\documentclass{amsart}
\usepackage{amsmath,amsthm}
\usepackage{amsfonts,amssymb}

\usepackage{enumerate}
  
\usepackage{hyperref} 
\newtheorem{theorem}{Theorem}
\newtheorem{Prop}[theorem]{Proposition}

\numberwithin{theorem}{section}

\theoremstyle{definition}
\newtheorem{Def}[theorem]{Definition}

\theoremstyle{remark}
\newtheorem{Opm}[theorem]{Remark}

\usepackage{appendix} 
 
\usepackage{float}

\usepackage{ dsfont }
\usepackage{ mathrsfs }
\usepackage{mathtools}

\def\Z{\mathds{Z} }

\def\C{\mathds{C} }

\def\K{\mathds{K} }

\DeclareMathOperator{\TKK}{TKK}
\DeclareMathOperator{\lTKK}{\overline{TKK}}
\DeclareMathOperator{\tTKK}{\widetilde{TKK}}
\DeclareMathOperator{\ad}{ad}
\DeclareMathOperator{\Ad}{Ad}
\DeclareMathOperator{\Inn}{Inn}
\DeclareMathOperator{\End}{End}
\DeclareMathOperator{\bessel}{\mathcal B_\lambda}

\DeclareMathOperator{\SB}{SB}

\DeclareMathOperator{\Lie}{Lie}
\newcommand{\iSB}{\SB^{-1}}

\newcommand{\bfip}[1]{\left<{#1}\right>_\mathcal B}

\newcommand{\bfipx}[1]{\left<{#1}\right>_{\mathcal B(x)}}
\newcommand{\ip}[1]{\left<{#1}\right>_W}

\newcommand{\bbf}[1]{\left<{#1}\right>_\beta}
\newcommand{\pt}[1]{\partial_{#1}}
\newcommand{\mf}[1]{\mathfrak{#1}}
\newcommand{\ds}[1]{\mathds{#1}}
\newcommand{\mc}[1]{\mathcal{#1}}

\newcommand{\pushright}[1]{\ifmeasuring@#1\else\omit\hfill$\displaystyle#1$\fi\ignorespaces}
\newcommand{\ol}[1]{\overline{#1}}

\newcommand{\ot}[1]{\widetilde{#1}}

\newcommand{\id}{{\rm id}}

\newcommand{\minus}{\scalebox{0.9}{{\rm -}}}
\newcommand{\plus}{\scalebox{0.6}{{\rm+}}}

\DeclareMathOperator{\tr}{tr}
\DeclareMathOperator{\Der}{Der}
\DeclareMathOperator{\pil}{\pi_\lambda}
\DeclareMathOperator{\rol}{\rho_\lambda}

\DeclareMathOperator{\OSp}{OSp}
\DeclareMathOperator{\Sp}{Sp}
\DeclareMathOperator{\Mp}{Mp}

\DeclarePairedDelimiter\abs{\lvert}{\rvert}%
\DeclarePairedDelimiter\norm{\lVert}{\rVert}%
\makeatletter
\let\oldabs\abs
\def\abs{\@ifstar{\oldabs}{\oldabs*}}
\let\oldnorm\norm
\def\norm{\@ifstar{\oldnorm}{\oldnorm*}}
\makeatother

\begin{document}
\title[A general approach to constructing minimal representations]{A general approach to constructing minimal representations of Lie supergroups}

\author{Sigiswald Barbier}
\address{Department of Electronics and Information Systems \\Faculty of Engineering and Architecture\\Ghent University\\Krijgslaan 281, 9000 Gent\\ Belgium.}
\email{Sigiswald.Barbier@UGent.be}

\author{Sam Claerebout}
\address{Department of Electronics and Information Systems \\Faculty of Engineering and Architecture\\Ghent University\\Krijgslaan 281, 9000 Gent\\ Belgium.}
\email{Sam.Claerebout@UGent.be}


\begin{abstract}
In this paper we describe an approach to generalise minimal representations to the super setting for Lie superalgebras obtained from Jordan superalgebras using the TKK construction.    
This approach was used successfully to construct a Fock model, a Schr\"odinger model and intertwining Segal-Bargmann transform for the orthosymplectic Lie supergroup $\OSp(p,q|2n)$ and the exceptional Lie supergroup $\mathbb{D}(2,1;\alpha)$. 
We also describe some obstacles to use this approach for the periplectic and queer Lie superalgebras. 
\end{abstract}

\maketitle


\section*{Introduction}

For a regular Lie group there exists a correspondence between coadjoint orbits and irreducible unitary representations, which is often referred to as the orbit method. Within this framework, a minimal representation of a semisimple Lie group is an irreducible unitary representation which corresponds to the minimal coadjoint orbit. For the technical details and exact definitions we refer to \cite{GanSavin}. 
Minimal representations tend to have multiple interesting realisations, with their own benefits and downsides and have been studied in many different settings \cite{VergneRossi, DvorskySahi, KO3, KM, HKM}.

A prominent example is the Segal-Shale-Weil representation of the metaplectic group $\Mp_{\ds R}(2n)$, which is a double cover of the symplectic group $\Sp_{\ds R}(2n)$. This representation is also known under the name oscillator representation or metaplectic representation.
 For the Segal-Shale-Weil representation we have the Schr\"odinger model (also called the $L^2$-model) and the Fock model, which are connected to each other by an intertwining integral operator called the Segal-Bargmann transform. For more information about these models of the Segal-Shale-Weil representation and the (classical) Segal-Bargmann transform we refer to \cite{Folland}. Analogous realisations can be constructed for a large group of other minimal representations of Lie groups and are then also called Schr\"odinger and Fock models. For simple Lie groups a unified construction of Schr\"odinger models was given in \cite{HKM}. Later, for hermitian Lie groups of tube type, a unified construction of Fock models was given in \cite{HKMO}. There, they also construct a generalised Segal-Bargmann transform, which intertwines their Fock model with the Schr\"odinger model in \cite{HKM}.

Lie supergroups are a generalisation of Lie groups that allow us to mathematically describe supersymmetry. In \cite{BC1} and \cite{BF} the question was asked whether the minimal representations can also be defined and constructed for Lie supergroups. In \cite{BF}, an answer was given for the orthosymplectic Lie supergroup $\OSp(p,q|2n)$. Essentially, a Schr\"odinger model of a minimal representation analogue was constructed using a generalisation of the approach in \cite{HKM}. Following the approach in \cite{HKMO}, a Fock model and intertwining Segal-Bargmann transform analogue was also constructed for $\OSp(p,2|2n)$ in \cite{BCD}.

One of the reasons the orthosymplectic supergroup was chosen as a starting point is that it generalises the indefinite orthogonal group ${\rm O}(p,q)$. The minimal representation of ${\rm O}(p,q)$ is well studied, see e.g. \cite{KO1, KO2, KO3, KM}. It is therefore easier to see how these results should extend to $\OSp(p,q|2n)$. 
In \cite{BC3} and \cite{BCARXIV}, a Schr\"odinger and a Fock model of a minimal representation and an intertwining Segal-Bargmann transform  were constructed for the exceptional Lie supergroup $\ds D(2,1;\alpha)$. Unlike the orthosymplectic supergroup, it has no clear classical counterpart.

For Lie groups, minimal representations are always unitary representations. However, the minimal representations constructed for $\OSp(p,q|2n)$ and $\ds D(2,1;\alpha)$ are not unitary. Moreover, there does not exist any unitary representation of $\OSp(p,q|2n)$ for nonzero $p$, $q$ and $n$, see \cite[Theorem 6.2.1]{NeebSalmasian}. These minimal representations of Lie supergroups are expected to be examples of ``superunitary'' representations for some new meaning of the word ``superunitary''. There already exists some definitions of superunitary representations which generalise the notion of unitary representations, see e.g. \cite{dGM} and \cite{Tuynman}, but it is still unclear if these are the correct definitions to use.

One of the main motivations of this paper is the idea that studying minimal representations of Lie superalgebras gives us a better understanding of what the definition of a superunitary representation should be.  In \cite{BCARXIV} it was proven that the minimal representation of $\ds D(2,1;\alpha)$ is, according to the definitions in \cite{dGM}, a superunitary representation. In contrast, the minimal representation of $\OSp(p,q|2n)$ does not seem to be superunitary for any existing definition.

The classical cases given in \cite{HKM} and \cite{HKMO}, together with the super cases given in \cite{BF}, \cite{BCD}, \cite{BC3} and \cite{BCARXIV}, give us a general idea of what these minimal representation ought to be for other Lie supergroups. In this paper we try to distil a general framework out of these cases that also works for other Lie supergroups. In particular, for a Lie supergroup we formulate an approach to construct:
\begin{itemize}
\item a minimal representation,
\item a Schr\"odinger model,
\item a Fock model,
\item and an intertwining Segal-Bargmann transform.
\end{itemize}
We also construct a ``super-inner product'' called the Bessel-Fischer product. In many ways, these super-inner products on supervector spaces replace the role of inner products on vector spaces. In both \cite{dGM} and \cite{Tuynman}, they play an important part in defining superunitary representations.

The Lie algebras in \cite{HKM} and \cite{HKMO} are seen as a TKK construction of a Jordan algebra. Therefore, we restrict ourselves to Lie supergroups where the corresponding Lie superalgebra is obtained from a unital Jordan superalgebra using the TKK construction. Such Lie superalgebras are then called TKK algebras. A complete list of these TKK algebras over the complex field, together with the corresponding Jordan superalgebras, is given in Table \ref{TabTKK}. Multiple TKK algebras can correspond to the same Jordan superalgebra $J$, which is why we have two columns, labelled $\lTKK(J)$ and $\tTKK(J)$. A brief introduction to these TKK constructions is given in Section \ref{SS_TKK} while a full overview can be found in \cite{BC2}. Note that we will not prove that our approach works for all of the Lie superalgebras in Table \ref{TabTKK}. On the contrary, in Sections \ref{SecPeri} and \ref{SecQueer}, we will prove that at least some part of our approach does not work for the periplectic and queer Lie superalgebras. Specifically, in these cases we can not use the Bessel-Fischer product as a super-inner product to obtain a superunitary minimal representation.

\begin{table}[H]
\centering\renewcommand*{\arraystretch}{1.4}
\caption{The TKK algebras. When $\lTKK(J)=\tTKK(J)$ we only wrote it once.}\label{TabTKK}
\begin{tabular}{|c|c|c|c|}
\hline
$J$ & $\lTKK(J)$ & $\tTKK(J)$ & Restrictions \\\hline\hline
$JGL(m|n)$ & $\mf{sl}(2m|2n)$ &  & $m\neq n$\\\hline
$JGL(m|m)$ & $\mf{psl}(2m|2m)$ & $\mf{pgl}(2m|2m)$ & $m>1$\\\hline
$JOSP(m|2n)$ & $\mf{osp}(4n|2m)$ &  & $(m,n)\neq(0,1)$\\\hline
$JSpin(m-3|2n)$ & $\mf{osp}(m|2n)$ &  & $m\geq 3$, $(m,n)\neq(4,0)$\\\hline
$JPe(n)$ & $\mf{spe}(2n)$ & $\mf{pe}(2n)$ & $n>1$\\\hline
$JQ(n)$ & $\mf{psq}(2n)$ & $\mf{pq}(2n)$ & $n>1$\\\hline
$D_\alpha$ & $D(2,1;\alpha)$ &  & $\alpha\neq\{0,-1\}$\\\hline
$D_{-1}$ & $\mf{psl}(2|2)$ & $D(2,1;-1)$ & \\\hline
$E$ & $E(7)$ &  & \\\hline
$F$ & $F(4)$ & & \\\hline
$JP(0,n-3)$ & $H(n)$ & $\ds K C\ltimes\ot H(n)$ & $n\geq 5$\\\hline
\end{tabular}
\end{table}

\subsection*{Structure of the paper}

The paper is organised as follows. 

In Section \ref{SecConstr}, we give the basic definitions, conventions and preliminary results. Specifically, in Sections \ref{SSgl}, \ref{SSpe} and \ref{SSq} we define the general linear, periplectic and queer superalgebras, respectively. In Section \ref{SS_Group}, we define Lie supergroups. Section \ref{SS_TKK} contains a brief introduction to TKK constructions based on the results in \cite{BC2}. In Section \ref{Section Special cases}, we specify the TKK algebras where our approach has already been applied as mentioned in the introduction. Sections \ref{Section Bessel operators} and \ref{SS_pol_real} are based on the results in \cite{BC1}. There, we introduce the Bessel operators (Section \ref{Section Bessel operators}) and a polynomial realisation from which our minimal representation can be constructed (Section \ref{SS_pol_real}).

Section \ref{SecMin} is the main body of this paper. In Section \ref{SS_Min_Rep_Def}, we give our definition of a minimal representation on the Lie superalgebra level. We then discuss a method to integrate our minimal representation to the group level in Section \ref{SS_Int}. Based on this discussion, we define a Schr\"odinger and Fock model in Sections \ref{Section Schrod model} and \ref{Section Fock model}, respectively. In Section \ref{SS_BF}, we define the Bessel-Fischer product and discuss some of the properties we expect it to have. In Section \ref{SSintop}, we construct an intertwining operator between Schr\"odinger and Fock models, which is then used in Section \ref{SS_SB} to define a Segal-Bargmann transform.

Lastly, in Sections \ref{SecPeri} and \ref{SecQueer}, we discuss our approach and shortcomings thereof for the periplectic and queer Lie superalgebras, respectively. Both of these sections are split into three parts. In the first part, Sections \ref{SS_Pe_TKK} and \ref{SS_Q_TKK}, we give the TKK construction explicitly. The second part, Sections \ref{SS_Pe_Lambda} and Sections \ref{SS_Q_Lambda}, discusses the character on which our representations depend. In the last part, Sections \ref{SS_Pe_Bessel} and Sections \ref{SS_Q_Bessel}, we give the Bessel operators explicitly and show that the Bessel-Fischer product does not have the desired properties.

\section{Preliminaries}\label{SecConstr}

In this paper we denote the complex unit by $\imath$ and use the convention that $\ds N = \{0,1,2,\ldots\}$. The field we work over can be either the field of real numbers $\ds R$ or the field of complex numbers $\ds C$ and will be denoted by $\ds K$. Function spaces will always be defined over $\ds C$.

A \textit{supervector space} is defined as a $\Z/2\Z$-graded vector space. This means that a supervector space $V$ has a decomposition $V=V_{\ol 0}\oplus V_{\ol 1}$. We call $V_{\ol 0}$ and $V_{\ol 1}$ the even part and the odd part of $V$, respectively.

An element $v\in V$ is called \textit{homogeneous} if $v\in V_i$, $i\in \Z/2\Z$. The \textit{parity} of a homogeneous element $v\in V$, denoted by $|v|$, is the $i\in \Z/2\Z$ such that $v\in V_i$. When we use $|v|$ in a formula, we are considering homogeneous elements, with the implicit convention that the formula has to be extended linearly for arbitrary elements.

We define the \textit{dimension} of a supervector space $V$ as $\dim(V) := (d_{\ol 0}|d_{\ol 1})$, where $d_i = \dim(V_i)$ are the dimensions of $V_i$ as vector spaces.

The supervector space $V$ with $V_{\ol 0} = \ds K^m$ and $V_{\ol 1} = \ds K^n$ is denoted by $\ds K^{m|n}$.

A \textit{superalgebra} is a supervector space $A=A_{\ol 0}\oplus A_{\ol 1}$ for which $A$ is an algebra and $A_iA_j\subseteq A_{i+j}$.

Let us now introduce the formal definition of a Lie superalgebra and a Jordan superalgebra.
\begin{Def}
A \textit{Lie superalgebra} is a superalgebra $\mf{g}=\mf{g}_{\ol 0}\oplus\mf{g}_{\ol 1}$ with a bilinear multiplication $[\cdot \,, \cdot]$ satisfying the following properties:
\begin{itemize}
\item Skew-supersymmetry: $[a,b]=-(-1)^{|a||b|}[b,a]$,
\item Super Jacobi identity: $[a,[b,c]]=[[a,b],c]+(-1)^{|a||b|}[b,[a,c]]$.
\end{itemize}
The multiplication on $\mf g$ is called its \textit{Lie bracket}.
\end{Def}

\begin{Def}
A \textit{Jordan superalgebra} is a supercommutative superalgebra $J = J_{\ol 0}\oplus J_{\ol 1}$ satisfying the Jordan identity
\begin{align*}
(-1)^{|x||z|}[L_{x}, L_{yz}]+(-1)^{|y||x|}[L_{y}, L_{zx}]+(-1)^{|z||y|}[L_{z}, L_{xy}]=0 \text{ for all } x,y,z \in J.
\end{align*}
Here the operator $L_x$ is left multiplication with $x$ and $[\cdot\,,\cdot]$ is the supercommutator, i.e.\ $[L_x,L_y] := L_xL_y - (-1)^{|x||y|}L_yL_x$. The algebra product on $J$ is called the \textit{Jordan product}.
\end{Def}

The following three subsections give some examples of Jordan and Lie superalgebras which we will need later on, particularly in Sections \ref{SecPeri} and \ref{SecQueer}. Other examples and more information can be found in e.g.\ \cite{CW}.

\subsection{The general linear superalgebras}\label{SSgl}
The endomorphisms of the supervector space $\ds K^{m|n}$ form an associative superalgebra. We can express it using matrices as 
\begin{align*}
\End(\ds K^{m|n}):= \left\lbrace \left(\begin{array}{c|c}
 a & b\\ \hline
 c & d\\
\end{array}\right)\mid a\in \ds K^{m\times m}, b\in \ds K^{m\times n}, c\in \ds K^{n\times m}, d\in \ds K^{n\times n} \right\rbrace,
\end{align*}
where the block diagonal matrices $a$ and $d$ give the even part, while the odd part is given by the off-diagonal blocks $b$ and $c$. 

As supervector spaces, the \textit{general linear Lie superalgebra} $\mf{gl}(m|n)$ and the \textit{full linear Jordan superalgebra} $JGL(m|n)$ are equal to $\End(\ds K^{m|n})$. 
The Lie bracket on $\mf{gl}(m|n)$ is given by \[ [x,y] := xy - (-1)^{|x||y|}yx \text{ for }x,y\in \mf{gl}(m|n),\]
 while  the Jordan product on $JGL(m|n)$ is given by 
 \[ x\cdot y := \frac{1}{2}(xy+(-1)^{|x||y|}yx) \text{ for } x,y\in JGL(m|n).\]

Let $E_{ij}$ be the $(m|n) \times (m|n)$-matrix where the only non-zero entry is a $1$ on the $i$th row and $j$th column. Then $\{E_{ij}: 1\leq i,j \leq m+n\}$ is a basis of $\mf{gl}(m|n)$ and $JGL(m|n)$.

\subsection{The periplectic superalgebras}\label{SSpe}

Let $\bbf{\cdot\, , \cdot}$ be a non-degenerate, supersymmetric, odd bilinear form on $\ds K^{n|n}$. In this paper we will always assume the corresponding matrix form is given by
\begin{align*}
\beta := \left(\begin{array}{c|c}
 0 & I_n\\ \hline
 I_n & 0\\
\end{array}\right),
\end{align*}
with $I_n$ the $n$-dimensional identity matrix.

The \textit{periplectic Lie superalgebra} $\mf {pe}(n)$ is the subalgebra of $\mf {gl}(n|n)$ consisting of all operators that leave $\beta$ invariant. We can express it using matrices as
\begin{align*}
\mf {pe}(n) := \left\lbrace \left(\begin{array}{c|c}
 a & b\\ \hline
 c & -a^t\\
\end{array}\right)\mid a, b, c\in \ds K^{n\times n}, b^t = b, c^t= -c\right\rbrace \subset \mf {gl}(n|n).
\end{align*}

The \textit{special periplectic Lie superalgebra} is defined as
\begin{align*}
\mf {spe}(n) := \left\lbrace \left(\begin{array}{c|c}
 a & b\\ \hline
 c & -a^t\\
\end{array}\right) \mid  a, b, c\in \ds K^{n\times n}, b^t = b, c^t= -c, \tr(a)=0\right\rbrace \subset \mf{pe}(n),
\end{align*}
which is simple for $n\geq 3$.

The \textit{Jordan-periplectic superalgebra} $JPe(n)$ is the subalgebra of $JGL(n|n)$ consisting of self-adjoint operators with respect to $\beta$. We can express it using matrices as
\begin{align*}
JPe(n) := \left\lbrace \left(\begin{array}{c|c}
 a & b\\ \hline
 c & a^t\\
\end{array}\right)\mid a, b, c\in \ds K^{n\times n}, b^t = -b, c^t=c \right\rbrace \subset JGL(n|n).
\end{align*}
Note that $b$ is symmetric and $c$ anti-symmetric for $\mf {pe}(n)$  and the other way around for $JPe(n)$.

\subsection{The queer superalgebras}\label{SSq}

The \textit{queer Lie superalgebra} $\mf {q}(n)$ is the subalgebra of $\mf gl(n|n)$ consisting of all operators in $\mf {gl}(n|n)$ commuting with the operator $\beta$. We can express it using matrices as
\begin{align*}
\mf q(n) := \left\lbrace \left(\begin{array}{c|c}
 a & b\\ \hline
 b & a\\
\end{array}\right)\mid a, b\in \ds K^{n\times n} \right\rbrace \subset \mf {gl}(n|n).
\end{align*}

The \textit{special queer Lie superalgebra} is defined as
\begin{align*}
\mf {sq}(n) := \left\lbrace \left(\begin{array}{c|c}
 a & b\\ \hline
 b & a\\
\end{array}\right)\mid a, b\in \ds K^{n\times n}, \tr(b)=0 \right\rbrace \subset \mf q(n).
\end{align*}

The \textit{projective queer Lie superalgebra} is defined as
\begin{align*}
\mf {pq}(n) := \mf{q}(n)/\left<I_{2n}\right>.
\end{align*}

The \textit{projective special queer Lie superalgebra} is defined as
\begin{align*}
\mf {psq}(n) := \mf{sq}(n)/\left<I_{2n}\right>,
\end{align*}
which is simple for $n \geq 3$.

We define the \textit{Jordan-queer superalgebra} $JQ(n)$ as the subalgebra of all operators in $JGL(n|n)$ commuting with the operator $\beta$. Hence as supervector spaces $JQ(n)= \mf q(n)$.

\subsection{Lie supergroups}\label{SS_Group}

In the existing literature there exists two equivalent ways to characterise Lie supergroups, one based on supermanifolds and one based on pairs. We will use the characterisation based on pairs and refer to \cite[Chapter 7]{CCF} for more information about these structures and how they are connected.

\begin{Def}
A \textit{Lie supergroup} $G$ is a pair $(G_0, \mf g)$ together with a morphism $\sigma \colon G_0\rightarrow \End(\mf g)$ where $G_0$ is a Lie group and $\mf g$ is a Lie superalgebra for which
\begin{itemize}
\item $\mf g_{\ol 0}$ is the Lie algebra of $G_0$, i.e.\ $\Lie(G_0)\cong \mf g_{\ol 0}$,
\item the morphism $\sigma$ extends the adjoint representation $\Ad$ of $G_0$ on $\mf g_{\ol 0}$, i.e.\ for all $g\in G_0$ we have $\left.\sigma(g)\right|_{\mf g_{\ol 0}} = \Ad(g)$,
\item for all $X\in \mf g_{\ol 0}$ and $Y\in \mf g$ we have
\begin{align*}
d\sigma(X)Y= \left.\dfrac{d}{dt}\sigma(\exp(tX))Y\right|_{t=0}=[X,Y].
\end{align*}
\end{itemize}
\end{Def}

\subsection{The TKK construction}\label{SS_TKK}

From now on we denote by $J$ a given unital Jordan superalgebra over the field $\ds K$ and suppose $\dim(J) = (m|n)$.

We can associate a $3$-graded Lie superalgebra to $J$ via the TKK construction. There are multiple TKK constructions in the literature, see \cite{BC2} for an overview, but for a unital Jordan superalgebra all TKK constructions are equivalent to the ones we will construct in this section. To define these TKK-constructions we must first define the following structures which are interesting in their own right.

The \textit{space of derivations} is defined as
\begin{align*}
\Der(J) := \left\lbrace D \in \mf{gl}(J) \mid D(xy) = D(x)y+(-1)^{|x||D|}x D(y), x,y \in J\right\rbrace
\end{align*}
and the \textit{space of inner derivations} is defined as
\begin{align*}
\Inn(J) := \mbox{span}_{\ds K}\left\lbrace [L_x, L_y]\mid x,y\in J\right\rbrace.
\end{align*}
Here the operator $L_x \in \mf{gl}(J)$ is left multiplication by $x$ and $[\cdot\,,\cdot]$ is the bracket of $\mf{gl}(J)$, i.e.\ $[L_x,L_y] = L_xL_y - (-1)^{|x||y|}L_yL_x$.

The \textit{structure algebra} is defined as
\begin{align*}
\mf{str}(J) = \{L_x \mid x\in J\} \oplus \Der(J)
\end{align*}
and the \textit{inner structure algebra} is defined as
\begin{align*}
\mf{istr}(J) = \{L_x \mid x\in J\} \oplus \Inn(J) = \mbox{span}_{\ds K}\left\lbrace L_x, [L_x, L_y]\mid x,y\in J\right\rbrace.
\end{align*}

Let $J^+$ and $J^-$ be two copies of $J$. As vector spaces we define the following two \textit{TKK algebras} of $J$:
\begin{align*}
\lTKK(J) &:= J^-\oplus \mf{istr}(J) \oplus J^+,\\
\tTKK(J) &:= J^-\oplus \mf{str}(J) \oplus J^+.\\
\end{align*}
The Lie bracket is defined as follows. We embed $\mf{istr}(J)$ or $\mf{str}(J)$ as a subalgebra of $\lTKK(J)$ or $\tTKK(J)$, respectively and set
\begin{align*}
[x,u] &= 2L_{xu}+2[L_x,L_u], &&& [x,y] &= [u,v] = 0,\\
[L_a,x] &= ax, &&& [L_a,u] &= -au,\\
[D, x] &= Dx, &&& [D, u] &= Du,
\end{align*}
for homogeneous $x,y\in J^+$, $u,v\in J^-$, $a \in J$ and $D\in \Inn(J)$ or $D\in \Der(J)$.

For $\ds K = \ds C$ a complete list unital finite dimensional simple Jordan superalgebras together with the corresponding structure algebras and TKK algebras is given in \cite[Section 6.1]{BC2}. Therefore, this list gives us precisely the (complex) Lie superalgebras for which our approach is applicable. For convenience, we restated it in the introduction, Table \ref{TabTKK} using our notational conventions. When $\lTKK(J)\cong \tTKK(J)$, we left the $\tTKK(J)$ box empty.

Most of the results in this paper are only stated for $\TKK(J) := \tTKK(J)$. The $\lTKK(J)$ case is analogous. From now on we denote by $\mf g = \TKK(J)$ the Lie superalgebra associated to $J$ via the TKK construction.

\subsection{Special cases}
\label{Section Special cases}
In this paper we will often distinguish the following cases.
\begin{itemize}
\item If $J$ is a simple Jordan algebra, then $\mf g$ is a Lie algebra and the resulting representations are the ones studied in \cite{HKM} and \cite{HKMO}. We will refer to this case as the classical case.
\item If $J$ is the spin factor Jordan superalgebra $JSpin(m-3|2n)$, then $\mf g$ is the orthosymplectic Lie superalgebra $\mf{osp}(m|2n)$ and the resulting representations are the ones studied in \cite{BF} and \cite{BCD}. We will refer to this case as the orthosymplectic case.
\item If $J$ is the Jordan superalgebra $D_\alpha$, then $\mf g$ is the exceptional Lie superalgebra $D(2,1;\alpha)$ and the resulting representations are the ones studied in \cite{BC3} and \cite{BCARXIV}. We will refer to this case as the $D(2,1;\alpha)$ case.
\end{itemize}
In Sections \ref{SecPeri} and \ref{SecQueer} we will also discuss the cases where $J$ is the Jordan-periplectic superalgebra $JPe(n)$ and the Jordan-queer superalgebra $JQ(n)$, respectively. We refer to these cases the periplectic and queer cases.

\subsection{The Bessel operators}
\label{Section Bessel operators}
In \cite{BC1} a polynomial realisation for $\mf g = \TKK(J)$ was constructed that is an analogue of the conformal representations considered in \cite{HKM}. In the orthosymplectic and $D(2,1;\alpha)$ cases this polynomial realisation was then used to define the minimal representation. The polynomial realisation depends on a character $\lambda \colon \mathfrak{str}(J) \to \ds K$ and a crucial role is played by the Bessel operators.

\begin{Def} Consider a character $\lambda \colon \mathfrak{str}(J) \to \ds K$. For any $u,v \in J^{\minus}$ we define $\lambda_{u} \in (J^{\plus})^\ast$ and $\widetilde{P}_{u,v} \in J^{\minus}\otimes (J^{\plus})^\ast $ by 
$$\lambda_{u}(x) := -2\lambda(L_{xu})$$
and
$$\widetilde{P}_{u,v}(x) := (-1)^{\abs{x}(\abs{u}+\abs{v})}(L_u L_v + (-1)^{|u||v|}L_vL_u-L_{uv})(x)$$
for all $x\in J^{\plus}$. Then we define the \textit{Bessel operator} as
\begin{align*}
\bessel = \sum_{i=1}^{m+n} \lambda_{z_{i}}\pt{z_{i}}+ \sum_{i,j=1}^{m+n} \widetilde{P}_{z_{i},z_{j}}\pt{z_{j}}\pt{z_{i}},
\end{align*}
\end{Def}
with $(z_i)_i$ a homogeneous basis of $J^-$.

From \cite[Proposition 4.2]{BC1} we obtain the following result.

\begin{Prop}\label{PropBesCom}
The family of operators $\bessel(x)$ for $x\in J^+$, supercommutes for fixed $\lambda$, i.e.
\begin{align*}
[\bessel(u), \bessel(v)] = 0,
\end{align*}
for $u,v\in J^+$.
\end{Prop}

\subsection{The polynomial realisation $\pil$}\label{SS_pol_real}

The \textit{space of superpolynomials} over $\ds K$ is defined as
\begin{gather*}
\mathcal P\big(\ds K^{m|n}\big):=\mathcal P \big(\ds K^{m}\big)\otimes_\C\Lambda\big(\ds K^{n}\big),
\end{gather*}
where $\mc P(\ds K^m)$ denotes the space of complex-valued polynomials over the field $\ds K$ in $m$ variables and $\Lambda(\ds K^{n})$ denotes the Grassmann algebra in $n$ variables. The variables of $\mc P(\ds K^m)$ and $\Lambda(\ds K^{n})$ are called \textit{even} and \textit{odd} variables, respectively.

Let $z = (z_i)_{i=1}^{m+n}$ denote the variables of $\mathcal P(\ds K^{m|n})$, then they satisfy the commutation relations
\begin{gather*}
z_iz_j = (-1)^{|z_i||z_j|}z_j z_i,
\end{gather*}
for $i,j\in\{1, \ldots, m+n\}$. We also define the \textit{space of homogeneous superpolynomials} of degree $k$ as
\begin{gather*}
\mc P_k\big(\K^{m|n}\big) := \big\{p\in\mc P\big(\K^{m|n}\big)\mid \ds E p = k p\big\},
\end{gather*}
where $\ds E := \sum_{i=1}^{m+n}z_i\pt{z_i}$ is the Euler operator.

From \cite[Section 4.1]{BC1} we obtain the following polynomial realisation $\pil$ of $\mf g = J^- \oplus  \mf{str}(J) \oplus J^+ $ on $\mc P(J)\cong \mc P(\ds K^{m|n})$. Let $(z_i)_i$ be a homogeneous basis of $J$ and $(z_i^\pm)_i$ the corresponding homogeneous bases of $J^\pm$. For $D\in \mf{str}(J)$ and $i\in\{1, \ldots, m+n\}$ we have
\begin{align*}
&\bullet\, \pil(z_i^-) = -2\imath z_i,\\
&\bullet\, \pil(D) = \lambda(D)+ \sum_{j=1}^{m+n} [D,z_j^-]\pt{z_j},\\
&\bullet\, \pil(z_i^+) = -\dfrac{\imath}{2}\bessel(z_i),
\end{align*}
where we identify the homogeneous bases $(z_i^\pm)_i$ with the variables $(z_i)_i$ of $\mc P(\ds K^{m|n})$.

\section{The minimal representation}\label{SecMin}

\subsection{Definition}\label{SS_Min_Rep_Def}

In the orthosymplectic and $D(2,1;\alpha)$ cases a representation was constructed that can be considered the analogue of a minimal representation in the classical case. These ``minimal representations'' were obtained as quotients of $\pil$ for specific values of $\lambda$. More specifically, in both cases a non-trivial subspace $V_\lambda$ of $\mc P_2(\ds K^{m|n})$ was found on which the Bessel operators act trivially and which is also a $\mf{str}(J)$-module. It then follows from the Poincaré-Birkhoff-Witt theorem that
\begin{align*}
\mc I_\lambda := U(J^-) V_\lambda=  \mc{P}(\ds \K^{m|n}) V_\lambda
\end{align*}
is a submodule of $\pi_\lambda$. Here $U(J^-)$ denotes the universal enveloping algebra of $J^-$. Following the methods in \cite{HKMO} for the classical case, the quotient representation of $\mf g$ on $\mc{P} (\mathds{K}^{m |n}) / \mc{I}_\lambda$ can be considered the analogue of the minimal representation.

In general, let $V_\lambda$ be a subspace of $\mc P(\ds K^{m|n})$ on which the Bessel operators act trivially and which is also a $\mf{str}(J)$-module. We can define the quotient representation of $\mf g$ on $\mc{P} (\mathds{K}^{m |n}) / \mc{I}_\lambda$ with
\begin{align*}
\mc I_\lambda := U(J^-) V_\lambda=  \mc{P}(\ds K^{m|n}) V_\lambda.
\end{align*}
For every choice of $V_\lambda$ there is now a corresponding representation of $\mf g$ on $\mc{P} (\mathds{K}^{m |n}) / \mc{I}_\lambda$, which we will still denote by $\pil$. Note that only for certain characters $\lambda$ we will have interesting choices for $V_\lambda$. A \textit{minimal representation} of $\mf g$ can then be defined as ``the smallest'' infinite dimensional representation we can obtain using this method, if it exists. Here we can make the definition of smallest representation formal by using the Gelfand--Kirillov dimension \cite[Section 7.3]{Mu}. 

In the cases we will look at in more detail in the rest of this paper we  only considered representations for which $V_\lambda  \subset \mc P_2(\ds K^{m|n})$, as these correspond to the minimal representations for the special cases mentioned in Section \ref{Section Special cases}. However, most of our approach is still viable for the other possible quotient representations obtained by considering other choices for $V_\lambda$. We hypothesise that these representations are no longer minimal representations but are similar to the representations studied in \cite{MS}. 

\subsection{Integration to the group level}\label{SS_Int}

In order to integrate a representation to the groups level we can make use of the framework of Harish-Chandra supermodules developed in \cite{Alldridge}.

\begin{Def}\cite[Definition 4.1]{Alldridge} Let $V$ be a complex supervector space and $G = (G_0, \mf g)$ a Lie supergroup. Denote by $K$ a maximal compact subgroup of $G_0$ and by $\mf k=\Lie(K)$ the Lie algebra of $K$. Then $V$ is a \textit{$(\mf g,K)$-module} if it consists of $\mf k$-finite vectors and the derived action of $K$ agrees with the $\mf k$-module structure:
\begin{align*}
d\pi_0(X)(v) = \left.\dfrac{d}{dt}\pi_0(\exp(tX))(v)\right|_{t=0} = d\pi(X)(v)\text{ for all } X\in\mf k,\, v\in V
\end{align*}
and
\begin{align*}
\pi_0(k)(d\pi(X)(v)) = d\pi(\Ad(k)(X))(\pi_0(k)(v)), \text{ for all } k\in K, X\in \mf g, v \in V,
\end{align*}
where $\pi_0$ is the $K$-representation and $d\pi$ the $\mf g$-representation. A $(\mf g,K)$-module is a \textit{Harish–Chandra supermodule} if it is finitely generated over $U(\mf g)$ and is $K$-multiplicity finite.
\end{Def}

In general, let $G = (G_0, \mf g)$ be a supergroup. We denote by $K$ a maximal compact subgroup of $G_0$ and by $\mf k$ the corresponding maximal Lie algebra of $\mf g_{\ol 0}$, i.e.\ $\mf k \cong\Lie(K)$. In \cite[Theorem 4.6]{Alldridge} it is shown that if $(\mf g,K)$ is a Harish–Chandra supermodule then it integrates to a unique smooth Fréchet representation of moderate growth for $G$. Therefore, to prove that a representation of $\mf g$ integrates to the group level, it is sufficient to show that it corresponds to a Harish–Chandra supermodule.

One of the requirement for the polynomial realisation $\pil$ to correspond with a Harish–Chandra supermodule is that it must consist of $\mf k$-finite vectors. In the following two subsections we give two methods which turn $\pil$ into a representation consisting of $\mf k$-finite vectors. The first method is based on the observation that $\pil$ is still well defined when we extend the action to smooth superfunctions. This results in a generalisation of the classical Schr\"odinger representation as studied in \cite{HKM}. The second method is based on the observation that $\pil$ consists of $\mf{str}(J)$-finite vectors. This results in a generalisation of the classical Fock representation as studied in \cite{HKMO}.

\subsection{The Schr\"odinger model}
\label{Section Schrod model}
Let us extend $\pil$  to smooth superfunctions and suppose we have found a smooth superfunction $\kappa$ which is also $\mf k$-finite vector. (Technically $\kappa$ should be a global section of the supermanifold we obtain by quotienting the supermanifold associated with $J$ by the ideal $\mathcal{I_\lambda}$, but this will not be relevant for the following.) Then we can define the representation
\begin{align*}
W_\lambda &:= U(\mf g)\kappa \mod \mc I_\lambda,
\end{align*}
where $U(\mf g)$ is the universal enveloping algebra of $\mf g$ and the $\mf g$-module structure is given by $\pil$. We call $W_\lambda$ a \textit{Schr\"odinger space} and this representation a \textit{Schr\"odinger model} of our representation. By \cite[Lemma 2.23]{HKM}, these spaces consists of $\mf k$-finite vectors. Note that while \cite{HKM} only considers the classical case, the proof of \cite[Lemma 2.23]{HKM} still works when replacing the Lie algebras with Lie superalgebras.

In \cite{HKM}, $\kappa(z)$ was chosen to be $\ot K_{\frac{\nu}{2}}(|z|)$ for the classical case. (Unless $\mf g \cong \mf{so}(p,q)$ with $p+q$ odd and $p,q\geq 3$, see Remark \ref{Remark_opq}.) Here $\ot K$ is a renormalised K-Bessel function, $|z|$ is the norm on $J$ induced by the trace form and $\nu$ is a structure constant of $J$. In particular, if $J$ is euclidean then $\ot K_{\frac{\nu}{2}}(|z|) = \exp(-\tr(z))$, where $\tr(z)$ denotes the trace function of $J$. In \cite{BF}, $\kappa(z)$ was also chosen to be $\ot K_{\frac{\nu}{2}}(|z|)$ for the orthosymplectic case. Here $|z|$ now denotes a superfunction that generalises the classical norm. In \cite{BC3}, $\kappa(z)$ was chosen to be $\exp(-\tr(z))$ for the $D(2,1;\alpha)$ case.
In Section \ref{SSintop} we will make a choice of $\kappa$ such that it will correspond to $1$ in the Fock model. 

\begin{Opm}\label{Remark_opq}
In the classical case $\ot K_{\frac{\nu}{2}}(|z|)$ is not $\mf k$-finite if $\mf g\cong \mf{so}(p,q)$ with $p+q$ odd and $p,q\geq 3$ (\cite[Theorem A]{HKM}). Moreover, if $p,q\geq 4$ and $p+q$ odd then there is no minimal representation for any group with Lie algebra $g \cong \mf{so}(p,q)$ (\cite[Theorem 2.13]{VD}). Similarly, in the orthosymplectic case $\ot K_{\frac{\nu}{2}}(|z|)$ is also not $\mf k$-finite for $\mf g\cong \mf{osp}(p,q|2n)$ with $\nu\not\in -2\ds N$, $p+q$ odd and $p,q\geq 3$ (\cite[Theorem 5.3]{BF}).
\end{Opm}

\subsection{The Fock model}
\label{Section Fock model}

Note that the action of $\mf{str}(J)$ on elements in $\mc P_{k}(\ds K^{m|n})$ returns elements in $\mc P_{k}(\ds K^{m|n})$, for all $k\in \ds N$. Since $\mc P_{k}(\ds K^{m|n})$ is finite dimensional, this implies that the action of $\mf g$ on $\mc P(\ds K^{m|n})$ consists of $\mf{str}(J)$-finite vectors. Therefore, if there exist an isomorphism of $\mf g$ that maps $\mf k$ into $\mf{str}(J)$, then we could twist $\pil$ in such a way that it consist of $\mf k$-finite vectors. Since finiteness is preserved under complexification, it is sufficient for there to exist an isomorphism of $\mf g_{\ds C}$ that maps $\mf k_{\ds C}$ into $\mf{str}(J_{\ds C})$. Here $\mf g_{\ds C}$, $\mf k_{\ds C}$ and $J_{\ds C}$ are the complexifications of $\mf g$, $\mf k$ and $J$, respectively. Suppose $\gamma\in \End(\mf g_{\ds C})$ is an isomorphism such that $\gamma(\mf k_{\ds C}) \subseteq \mf{str}(J_{\ds C})$. Then we can define the representation
\begin{align*}
F_\lambda &:= U(\mf g)1 \mod \mc I_\lambda,
\end{align*}
where the $\mf g$-module structure is given by $\rol := \pil\circ \gamma$. We call $F_\lambda$ a \textit{Fock space} and this representation a \textit{Fock model} of our representation.
Note that $F_\lambda$ not only depends on the character $\lambda$ but also on the choice of $V_\lambda$ (or thus $\mc I_\lambda$) and on the choice of the twist $\gamma$. 

In the classical case with $J$ euclidean, the orthosymplectic case for $q=2$ and the $D(2,1;\alpha)$ case we can choose $\gamma$ to be the Cayley transform $c \in \End(\mf g_\C)$ defined as
\begin{align*}
c := \exp\left(\frac{\imath}{2}\ad(e^-))\exp(\imath\ad(e^+)\right),
\end{align*}
where $e^\pm$ denotes the unit of $J^\pm$. The corresponding Fock spaces for the classical, orthosymplectic and $D(2,1;\alpha)$ cases can be found \cite{HKMO}, \cite{BCD} and \cite{BC3}, respectively. Note that in \cite{HKMO} they use different conventions for their $\TKK$-constructions than the ones described in this paper, which results in a slightly different looking Cayley transform.

In \cite{BCARXIV} the Fock model was integrated to the group level for the $D(2,1;\alpha)$ case. It was then proven to be a superunitary representation using the definition of ``superunitary'' found in \cite{dGM}.

\begin{Opm}\label{Remark_Cayley}
In the classical case, if $J$ is euclidean then the Cayley transform $c$ induces an isomorphism between $\mf k_{\ds C}$ and $\mf{str}(J_{\ds C})$, which is why only euclidean Jordan algebras are considered in \cite{HKMO}. Similarly, in the orthosymplectic case $c$ only maps $\mf k_{\ds C}$ into $\mf{str}(J_{\ds C})$ when $J_{\ol 0}$ is euclidean, which is why only $q=2$ is considered in \cite{BC3}. In the $D(2,1;\alpha)$ case $J_{\ol 0}$ is always euclidean. The Cayley transform $c$ can still be used for non-euclidean Jordan superalgebras to define a Fock-like representation $\rho_\lambda^c := \pil\circ c$. However, this representation is, in general, not unitary.
\end{Opm}

\subsection{The Bessel-Fischer product}\label{SS_BF}

A significant difference between the super and the non-super cases is the roll of inner products on Schr\"odinger and Fock spaces. On the Fock space, our approach will be to generalise the following inner product from the classical case. Suppose $p,q\in \mc P(\ds C^m)$ and let $p(\bessel)$ be the operator obtained by replacing every instance of $z_i$ by $\bessel (z_i)$ in $p$. The Bessel-Fischer inner product of $p$ and $q$ is then defined as
\begin{align*}
\bfip{p,q} := \left. p(\bessel)\overline{q(\bar z)}\right|_{z=0}.
\end{align*}
It was introduced in \cite[Section 2.3]{HKMO} as an inner product on the polynomial space $\mc P(\C^m)$. In the classical case, it was proven that the Bessel-Fischer inner product is equal to the $L^2$-inner product on the Fock space \cite[Proposition 2.6]{HKMO}. In the orthosymplectic and $D(2,1;\alpha)$ cases this product was used as the starting point to generalise the Fock space to superspace. The Bessel-Fischer product is no longer an inner product in this case. However, under the right conditions it is a non-degenerate, sesquilinear, superhermitian form on $F_\lambda$, which is consistent with the definition of Hilbert superspaces given in \cite{dGM}.

In general, assume we have Bessel operators acting on $\mathcal{P}(J) \cong \mathcal P(\ds K^{m|n})$ as in Section \ref{Section Bessel operators}. Recall that this action depends on a character $\lambda\colon \mathfrak{str} (J) \to \ds K$. Then we define the Bessel-Fischer product on a superpolynomial space as follows.
\begin{Def}
For $p, q\in \mathcal P(\ds K^{m|n})$ we define the \textit{Bessel-Fischer product} of $p$ and $q$ as
\begin{align*}
\bfip{p,q} := \left. p(\bessel)\bar q(z)\right|_{z=0},
\end{align*}
where $\bar q(z) = \overline{q(\bar z)}$ is obtained by conjugating the coefficients of the polynomial $q$ and $p(\bessel)$ is obtained by replacing the occurrences of $z_i$ with $\bessel(z_i)$ for all $i\in \{1, \ldots, m+n\}$.
\end{Def}

In the orthosymplectic and $D(2,1;\alpha)$ cases the Bessel-Fischer product is a non-degenerate superhermitian sesquilinear form when restricted to the Fock space $F_\lambda$. Moreover, the Fock representation $\rol$ is skew-supersymmetric with respect to the Bessel-Fischer product. We conjecture that these properties hold in a more general setting. However, these properties can not be expected to hold for all unital Jordan superalgebras. For example, in Sections \ref{SecPeri} and \ref{SecQueer} we show that the Bessel-Fischer product is degenerate when $J$ is either the periplectic or queer Jordan superalgebra.

\subsection{An intertwining operator}\label{SSintop}

In the cases where we work with the Cayley transform $c$ for the twist $\gamma$ considered in the Fock model, see Section \ref{Section Fock model}, we can define $C := \exp(\frac{\imath}{2}e^-)\exp(\imath e^+)$. Then the Cayley transform $c$ can be rewritten as
\begin{align*}
c &= \exp\left(\frac{\imath}{2}\ad(e^-)\right)\exp\left(\imath\ad(e^+)\right) = \Ad\left(\exp(\frac{\imath}{2}e^-)\right)\Ad\left(\exp(\imath e^+)\right) = \Ad(C).
\end{align*}
Therefore, we formally have,
\begin{align*}
\rol(X) &= \pil(c(X)) = \pil(\Ad(C)X) = \pil(C)\pil(X)\pil(C)^{-1},
\end{align*}
which shows that $\pil(C)$ intertwines the actions of $\pil$ and $\rol$.

In general, let $(G_0, \mf g)$ be a supergroup and suppose we have $\gamma = \Ad(\Gamma)$ for some $\Gamma \in G_0$, then $\pil(\Gamma)$ intertwines the actions of $\pil$ and $\rol = \pil\circ \gamma$. In particular, with every Fock model we can then associate a Schr\"odinger model by choosing the superfunction $\kappa$ from Section \ref{Section Schrod model} to be $\pil(\Gamma)^{-1}1$.

\subsection{The Segal-Bargmann Transform}\label{SS_SB}

For the classical case, an operator that intertwines the Schr\"odinger and Fock model was constructed in \cite[Section 3]{HKMO} called the Segal-Bargmann transform.  Let $x= (x_i)_{i=1}^{m+n}$ and $z = (z_i)_{i=1}^{m+n}$ denote the variables of $\mathcal P(\ds R^{m|n})$ and $\mathcal P(\ds C^{m|n})$, respectively. Similarly, let $e_x$ and $e_z$ denote the elements corresponding with the unit of $J$ in $\mc{P}(\ds R^{m|n})$ and $\mc{P}(\ds \C^{m|n})$, respectively. Then this Segal-Bargmann transform is given in our framework by the transform $\SB \colon W_\lambda \rightarrow F_\lambda$ defined by
\begin{align}\label{EqSB}
\SB(f(x))(z) := \exp(-e_z)\ip{f(x),\ds I_\lambda(2x,2z)\exp(-2e_x)}
\end{align}
and its inverse is given by
\begin{align*}
\iSB(p(z))(x) = \exp(-2e_x)\bfip{p(z),\ds I_\lambda(2z,2x)\exp(-e_z)}.
\end{align*}
Here $\ip{\cdot\, ,\cdot }$ denotes $L^2$ inner product on $W_\lambda$ and $\ds I_\lambda(z,x)$ is the reproducing kernel on the Fock space constructed in \cite[Section 2.4]{HKMO}. In particular, we have
\begin{align}\label{EqRepKer}
\bfip{p(z),\ds I_\lambda(z,x)} = p(x) \quad \text{ for all } p\in F_\lambda.
\end{align}
For the orthosymplectic case, an intertwining Segal-Bargmann transform was constructed in \cite[Section 6]{BCD} and it is also given by definition (\ref{EqSB}). The reproducing kernel $\ds I_\lambda(z,x)$ is now the superfunction constructed in \cite[Section 4.2]{BCD} and it also satisfies equation (\ref{EqRepKer}).

For the $D(2,1;\alpha)$ case the Segal-Bargmann transform could also be defined using equation (\ref{EqSB}). However, in \cite[Section 7]{BC3} an intertwining Segal-Bargmann transform was constructed in a more easily generalisable way. First, the Bessel-Fischer product was used in combination with $\pil(C)$ as in Section \ref{SSintop} to define a corresponding product on the Schr\"odinger space $W_\lambda$ as follows. We set
\begin{align*}
\ip{f,g} := \bfip{\pil(C)f,\pil(C)g},
\end{align*}
for all $f,g\in W_\lambda$. The Segal-Bargmann transform is now defined by
\begin{align*}
\SB (f(x))(z) &:= \bfipx{\pil(C)f(x), \ds I_\lambda(x,z)}\\
&= \ip{f(x), \pil(C)^{-1}\ds I_\lambda(x,z)}\\
&=\ip{f(x), \exp(-\dfrac{1}{2}\bessel(e_x))\exp(-e_x)\ds I_\lambda(x,z)},
\end{align*}
where $\bfipx{\cdot\, , \cdot}$ denotes the Bessel-Fischer product in the $x$ variable and $\ds I_\lambda(z,x)$ is the superfunction constructed in \cite[Section 5.2]{BC3} satisfying equation (\ref{EqRepKer}). The advantage this definition has over the classical definition is that the intertwining and superunitarity properties are immediate from the definitions of $\pil(C)$ and $\ip{\cdot\, ,\cdot}$. Since the product on $W_\lambda$ is defined using the Bessel-Fischer product on $F_\lambda$, it is more convenient to use the inverse Segal-Bargmann transform, which is given by
\begin{align*}
\iSB (p(z))(x) &= \pil(C)^{-1}(\bfip{p,\ds I_\lambda(z,x)})\\
&= \exp(-\dfrac{1}{2}\bessel(e_x))\exp(-e_x)\bfip{p,\ds I_\lambda(z,x)},
\end{align*}
for all $p\in F_\lambda$.

In general, suppose we have a Schr\"odinger representation $\pil$, a Fock representation $\rol = \pil\circ \gamma$ and an element $\Gamma \in G_0$ such that $\gamma = \Ad(\Gamma)$. Let $F_\lambda$ be the complex Fock space, $W_\lambda$ be the real Schr\"odinger space with $\kappa = \pil(\Gamma)^{-1}1$ and $\ds I_\lambda(\cdot\, ,\cdot )$ be a superfunction that satisfies equation (\ref{EqRepKer}). Then we can define the following
\begin{Def}\label{DefSBinv}
The \textit{inverse Segal-Bargmann transform} $\iSB: F_\lambda \rightarrow W_\lambda$ is defined by
\begin{align*}
\iSB (p(z))(x) &:= \pil(\Gamma)^{-1}(\bfip{p,\ds I_\lambda(z,x)}),
\end{align*}
for all $p\in F_\lambda$.
\end{Def}

Define a sesquilinear form $\ip{\cdot \, ,\cdot}$ on $W_\lambda$ by
\begin{align*}
\ip{f,g} := \bfip{\pil(\Gamma)f,\pil(\Gamma)g},
\end{align*}
for all $f,g\in W_\lambda$. Then $\iSB$ automatically preserves the bilinear forms and we can define the Segal-Bargmann transform as follows.

\begin{Def}\label{DefSB}
The \textit{Segal-Bargmann transform} $\SB\colon  W_\lambda \rightarrow F_\lambda$ is defined by
\begin{align*}
\SB (f(x))(z) &:= \ip{f(x), \pil(\Gamma)^{-1}\ds I_\lambda(x,z)}\\
\end{align*}
for all $f\in W_\lambda$.
\end{Def}

Note that $\SB\circ \iSB = \id_{F_\lambda}$ and $\iSB \circ \SB = \id_{W_\lambda}$ as desired.

\section{The periplectic case}\label{SecPeri}

We now give the TKK constructions and Bessel operators of Section \ref{SecConstr} explicitly for the specific case where $J$ is the Jordan-periplectic superalgebra $JPe(n)$, with $n\geq 2$. The purpose of this section is to show that our approach does not fully work in this case. Specifically, we will show that the Bessel-Fischer product can not be used as a non-degenerate form.

At the start of Section \ref{SecMin} we remarked that different choices for $V_\lambda$ are possible, which could lead to different infinite dimensional representations. However, as we will show, for the periplectic case the Bessel-Fischer product is either degenerate on the associated Fock space or the representation is finite dimensional. This holds regardless of the chosen $V_\lambda$.

\subsection{The TKK construction}\label{SS_Pe_TKK}

From the table in Section \ref{SS_TKK} we see that the $\TKK$-construction of $JPe(n)$ is isomorphic to the periplectic Lie superalgebra $\mf{pe}(2n)$. We will first show this statement more explicitly.

Recall the definitions of the superalgebras given in Section \ref{SSpe} and the elements $E_{ij}$ in Section \ref{SSgl}. A basis of the Jordan-periplectic superalgebra $JPe(n)$ is given by
\begin{align*}
x_{ij} &:= E_{ij}+E_{j+n,i+n}, &&\mbox{ for } &&1\leq i,j\leq n,\\
\beta_{ij} &:= E_{i,j+n}-E_{j,i+n}, &&\mbox{ for } &&1\leq i < j \leq n,\\
\gamma_{ij} &:= E_{i+n, j}+E_{j+n,i}, &&\mbox{ for } &&1\leq i \leq j \leq n. 
\end{align*}
The Jordan multiplication in terms of the basis elements is given by
\begin{align*}
x_{ij}\cdot x_{kl} &= \dfrac{1}{2}(\delta_{jk}x_{il}+\delta_{il}x_{kj}), & x_{ij}\cdot \beta_{kl} &= \dfrac{1}{2}(\delta_{jk}\beta_{il}-\delta_{jl}\beta_{ik}),\\
x_{ij}\cdot \gamma_{kl} &= \dfrac{1}{2}(\delta_{ik}\gamma_{jl}+\delta_{il}\gamma_{jk}), & \beta_{ij}\cdot \gamma_{kl} &= \dfrac{1}{2}(\delta_{jk}x_{il}+\delta_{jl}x_{ik}-\delta_{ik}x_{jl}-\delta_{il}x_{jk}),\\
\beta_{ij}\cdot \beta_{kl} &= \gamma_{ij}\cdot \gamma_{kl} = 0.
\end{align*}

A basis of the periplectic Lie superalgebra $\mf g = \mf{pe}(2n)$ is given by
\begin{align*}
\ol x_{ij} &:= E_{ij}-E_{j+2n,i+2n}, &&\mbox{ for } &&1\leq i, j \leq 2n,\\
\ol \beta_{ij} &:= E_{i,j+2n}+E_{j,i+2n}, &&\mbox{ for } &&1\leq i < j \leq 2n,\\
\ol \gamma_{ij} &:= E_{i+2n, j}-E_{j+2n,i}, &&\mbox{ for } &&1\leq i \leq j \leq 2n. 
\end{align*}
The Lie bracket in terms of the basis elements is given by
\begin{align*}
[\ol x_{ij}, \ol x_{kl}] &= \delta_{jk}\ol x_{il}-\delta_{il}\ol x_{kj}, & [\ol x_{ij}, \ol \beta_{kl}] &= \delta_{jk}\ol \beta_{il}+\delta_{jl}\ol \beta_{ik},\\
[\ol x_{ij},\ol \gamma_{kl}] &= \delta_{ik}\ol \gamma_{lj}-\delta_{il}\ol \gamma_{kj}, & [\ol \beta_{ij}, \ol \gamma_{kl}] &=\delta_{jk}\ol x_{il}-\delta_{jl}\ol x_{ik}+\delta_{ik}\ol x_{jl}-\delta_{il}\ol x_{jk},\\
[\ol\beta_{ij},\ol \beta_{kl}] &= [\ol\gamma_{ij},\ol \gamma_{kl}] = 0.
\end{align*}
We have the $3$-grading $\mf g = \mf g_- \oplus \mf g_0 \oplus \mf g_+$ with
\begin{align*}
\mf g_- &= \mbox{span}_{\ds K}\left\lbrace\ol x_{2i-1,2j}, \ol \beta_{2i-1,2j-1}, \ol \gamma_{2i,2j}\mid 1\leq i,j\leq n\right\rbrace,\\
\mf g_+ &= \mbox{span}_{\ds K}\left\lbrace\ol x_{2i,2j-1}, \ol \beta_{2i,2j}, \ol \gamma_{2i-1,2j-1}\mid1\leq i,j\leq n\right\rbrace,\\
\mf g_0 &= \mbox{span}_{\ds K}\left\lbrace\ol x_{2i,2j}, \ol x_{2i-1,2j-1}, \ol \beta_{2i,2j-1}, \ol \gamma_{2i,2j-1}\mid1\leq i,j\leq n\right\rbrace.
\end{align*}

An explicit isomorphism $\phi$ between $ JPe(n)^-\oplus \mf{istr}(JPe(n))\oplus JPe(n)^+$ 	and $\mathfrak{spe} (2n) $ is uniquely determined by
\begin{align*}
& \phi(x_{ij}^-) = \ol x_{2j-1,2i}, &  & \phi(x_{ij}^+) = \ol x_{2j,2i-1}, &  &\phi(2L_{x_{ij}}) = \ol x_{2j,2i}-\ol x_{2j-1,2i-1},& \\
&  \phi(\beta_{ij}^-)= -\ol \gamma_{2i,2j} , & &\phi(\beta_{ij}^+) = \ol \gamma_{2i-1,2j-1}, &  &\phi(2L_{\beta_{ij}}) = \ol \gamma_{2i-1,2j}+\ol \gamma_{2i,2j-1},\\
& \phi(\gamma_{ij}^-) = \ol \beta_{2i-1,2j-1}, &  &\phi(\gamma_{ij}^+) = -\ol \beta_{2i,2j}, &  &\phi(2L_{\gamma_{ij}}) = \ol \beta_{2i,2j-1}+\ol \beta_{2i-1,2j}.
\end{align*}
We define the derivative $D_p\in \mf{str}(JPe(n))$ by
\begin{align*}
&D_p(x_{ij}) := 0, && D_p(\beta_{ij}) := -2\beta_{ij}, && D_p(\gamma_{ij}) := 2\gamma_{ij},
\end{align*}
for all $i,j$. Then, we have $\mf{str}(JPe(n)) = \{D_p\} \oplus \mf{istr}(JPe(n))$ and $ \phi(D_p) =\sum_i \ol x_{ii}$ extends the isomorphism $\phi$ to an isomorphism between $\TKK (J)$ and $\mf g =\mathfrak{pe}(2n)$. 

If we use the standard basis $\{E_{ij}: 1\leq i,j \leq 2n\}$ for $\mf{gl}(n|n)$, then $\mf g_0\cong \mf{gl}(n|n)$ by the isomorphism
\begin{align*}
\ol x_{2i,2j}\mapsto E_{i,j} && \ol x_{2i-1,2j-1}\mapsto -E_{j+n,i+n} && \ol \beta_{2i,2j-1} \mapsto E_{i,j+n}
\end{align*}
for $1\leq i,j \leq n$ and
\begin{align*}
\ol \gamma_{2i,2j-1}\mapsto -E_{j+n,i}\quad \text{ for } i\geq j, && \ol \gamma_{2i,2j-1}\mapsto E_{i+n,j}\quad \text{ for } i < j.
\end{align*}

\subsection{The character $\lambda$}\label{SS_Pe_Lambda}

We now look at a character $\lambda$ on $\mf g_0$.  Since $\lambda$ is an even Lie superalgebra morphisms, we immediately have $\lambda(\ol \beta_{2i,2l-1}) = 0$ and $\lambda(\ol \gamma_{2l-1,2j})= 0$.
We have  $\lambda([a,b])=\lambda(a)\lambda (b) -(-1)^{\abs{a}\abs{b}} \lambda(b) \lambda (a)=0$, and thus
\begin{align*}
\lambda([\ol x_{2i,2j}, \ol x_{2k,2l}]) &= 0, & \lambda([\ol x_{2i-1,2j-1}, \ol x_{2k-1,2l-1}]) &= 0,\\
\lambda([\ol \beta_{2i,2j-1}, \ol \gamma_{2k,2l-1}]) &= 0.
\end{align*}
Therefore, we have
\begin{align*}
\delta_{kj}\lambda(\ol x_{2i,2l}) &= \delta_{il}\lambda(\ol x_{2k,2j}), & \delta_{kj}\lambda(\ol x_{2i-1,2l-1}) &= \delta_{il}\lambda(\ol x_{2k-1,2j-1}), \\
\delta_{jl}\lambda(\ol x_{2i,2k}) &= \delta_{ik}\lambda(\ol x_{2j-1,2l-1}).
\end{align*}

For all $1\leq i ,j\leq n$ we now conclude
\begin{align*}
\lambda(\ol x_{2i-1,2j-1}) &= \lambda(\ol x_{2i,2j}) =\delta_{ij} \lambda(\ol x_{11}),\\
\lambda(\ol \beta_{2i-1,2j}) &= \lambda(\ol \gamma_{2i-1,2j}) = 0.
\end{align*}
Using the isomorphism $\phi$ we now find for a character $\lambda$ on $\mf{str}(JPe(n))$ that
\begin{align*}
\lambda(L_{x_{ij}}) = \lambda(L_{\beta_{ij}}) = \lambda(L_{\gamma_{ij}}) = 0
\end{align*}
and $\lambda(D_p) = \lambda_p := 2n \lambda(\phi^{-1}(\ol x_{11}))$. In particular, the only character on $\mf{istr}(JPe(n))$ is the trivial character $\lambda=0$.

\subsection{The Bessel operators}\label{SS_Pe_Bessel}
A technical but straightforward calculation gives us
\begin{align*}
\bessel(x_{kl}) &= \sum_{i,j=1}^n (x_{ij}\pt{x_{lj}}+\beta_{ij}\pt{\beta_{lj}})\pt{x_{ik}} +(1+\delta_{ik})(\gamma_{ij}\pt{x_{lj}} + x_{ji}\pt{\beta_{lj}})\pt{\gamma_{ik}},\\
\bessel(\beta_{kl}) &=\sum_{i,j=1}^n x_{ij}(1+\delta_{jl})\pt{x_{ik}}\pt{\gamma_{jl}}-x_{ji}(1+\delta_{ik})\pt{\gamma_{ik}}\pt{x_{jl}}\\
&\quad +\gamma_{ij}(1+\delta_{ik})(1+\delta_{jl})\pt{\gamma_{ik}}\pt{\gamma_{jl}},\\
\bessel(\gamma_{kl}) &=\sum_{i,j=1}^n (\gamma_{ij}\pt{x_{ki}}+ x_{ij}\pt{\beta_{ki}})\pt{x_{lj}} +(x_{ji}\pt{x_{ki}} - 2\beta_{ij}\pt{\beta_{ki}})\pt{\beta_{lj}}.
\end{align*}

Note that because $\lambda$ is the trivial character on $\mf{istr}(JPe(n))$, we do not have the terms with only one derivative. In particular, all terms in the Bessel operators derive what it acts on twice. Therefore $\bessel(x)p = 0$, for all $x\in J$ and $p\in \mc P_{1}(J)$. This implies that every element of $\mc P_{1}(J)$ is a degenerate element for the Bessel-Fischer product on $\mc P(J)$. In order for the Bessel-Fischer to be non-degenerate on $F_\lambda$ we must therefore have $\mc P_{1}(J)\subseteq V_{\lambda}$, but then
\begin{align*}
\bigoplus_{k=1}^{\infty}\mc P_k(J) \subseteq \mc I_\lambda,
\end{align*}
which implies $F_\lambda$ is isomorphic to $\ds C$ or $0$. Either way, we no longer have an infinite dimensional representation.

Let $J$ now once again be an arbitrary unital Jordan superalgebra. The arguments above show that if $\lambda$ is the trivial character on $\mf{istr}(J)$ then either $F_\lambda$ is isomorphic to $\ds K$ or $0$ or the Bessel-Fischer product is degenerate on $F_\lambda$.

\section{The queer case}\label{SecQueer}

We will now look in detail at the case where $J$ is the Jordan-queer superalgebra $JQ(n)$, with $n\geq 2$. Analogous to Section \ref{SecPeri}, the purpose of this section is to show that the Bessel-Fischer product can not be used as a non-degenerate form.

\subsection{The TKK construction}\label{SS_Q_TKK}

From the table in Section \ref{SS_TKK} we see that the $\TKK$-construction of $JQ(n)$ is isomorphic to the projective queer Lie superalgebra $\mf{pq}(2n)$. In this section we show this statement more explicitly.

Recall the definitions of the superalgebras given in Section \ref{SSq} and the elements $E_{ij}$ in Section \ref{SSgl}. A basis of the Jordan-queer superalgebra $JQ(n)$ is given by
\begin{align*}
y_{ij} &:= E_{ij}+E_{i+n,j+n}, &&\mbox{ for } &&1\leq i,j\leq n,\\
\theta_{ij} &:= E_{i+n,j}+E_{i,j+n}, &&\mbox{ for } &&1\leq i,j\leq n.
\end{align*}
The Jordan multiplication in terms of the basis elements is given by
\begin{align*}
y_{ij}\cdot y_{kl} &= \dfrac{1}{2}(\delta_{jk}y_{il}+\delta_{il}y_{kj}), & y_{ij}\cdot \theta_{kl} &= \dfrac{1}{2}(\delta_{jk}\theta_{il}+\delta_{il}\theta_{kj}),\\
\theta_{ij}\cdot \theta_{kl} &= \dfrac{1}{2}(\delta_{jk}y_{il}-\delta_{il}y_{kj}).
\end{align*}

The projective queer Lie superalgebra is given by
\begin{align*}
\mf g = \mf{pq}(2n) = \mbox{span}_{\ds K}\left\lbrace y_{ij}, \theta_{ij} : 1\leq i,j\leq 2n\right\rbrace/\left<\sum y_{ii}\right>.
\end{align*}
The Lie bracket in terms of the basis elements is given by
\begin{align*}
[y_{ij}, y_{kl}] &= \delta_{jk}y_{il}-\delta_{il}y_{kj}, & [y_{ij}, \theta_{kl}] &= \delta_{jk}\theta_{il}-\delta_{il}\theta_{kj},\\
[\theta_{ij}, \theta_{kl}] &= \delta_{jk}y_{il}+\delta_{il}y_{kj}.
\end{align*}
We have the $3$-grading $\mf g = \mf g_- \oplus \mf g_0 \oplus \mf g_+$ with
\begin{align*}
\mf g_- &= \mbox{span}_{\ds K}\left\lbrace y_{2i-1,2j}, \theta_{2i-1,2j}\mid 1\leq i,j\leq n\right\rbrace,\\
\mf g_+ &= \mbox{span}_{\ds K}\left\lbrace y_{2i,2j-1}, \theta_{2i,2j-1}\mid1\leq i,j\leq n\right\rbrace,\\
\mf g_0 &= \mbox{span}_{\ds K}\left\lbrace y_{2i,2j}, \theta_{2i,2j}, y_{2i-1,2j-1}, \theta_{2i-1,2j-1}\mid1\leq i,j\leq n\right\rbrace/\left<\sum y_{ii}\right>.
\end{align*}

An explicit isomorphism $\phi$ between $\TKK(JQ(n)) = JQ(n)^-\oplus \mf{istr}(JQ(n))\oplus JQ(n)^+$ and $\mf {psq}(2n)$ is uniquely determined by
\begin{align*}
&\phi(y_{ij}^-) = y_{2i-1,2j}, && \phi(\theta_{ij}^-) = \theta_{2i-1,2j},\\
& \phi(y_{ij}^+) = y_{2i,2j-1}, && \phi(\theta_{ij}^+) = \theta_{2i,2j-1},\\
&\phi(2L_{y_{ij}}) = y_{2i,2j}-y_{2i-1,2j-1}, && \phi(2L_{\theta_{ij}}) = \theta_{2i,2j}-\theta_{2i-1,2j-1}.
\end{align*}
We define the derivative $D_q\in \mf{str}(JQ(n))$ by
\begin{align*}
&D_q(y_{ij}) := 0, && D_q(\theta_{ij}) = 2y_{ij},
\end{align*}
for all $i,j$. Then, we have $\mf{str}(JQ(n))=\{D_q\}\oplus \mf{istr}(JQ(n))$ and $\phi(D_q) = \sum_i \theta_{ii}$ extends the isomorphism $\phi$ to an isomorphism between $\TKK(J)$ and $\mf g = \mf{pq}(2n)$.

\subsection{The character $\lambda$}\label{SS_Q_Lambda}

We now look at a character $\lambda$ on $\mf g_0$. From 
\begin{align*}
\lambda([y_{2i,2j}, y_{2k,2l}]) &= 0, & \lambda([y_{2i-1,2j-1}, y_{2k-1,2l-1}]) &= 0,\\
\lambda([\theta_{2i,2j}, \theta_{2k,2l}]) &= 0, & \lambda([\theta_{2i-1,2j-1}, \theta_{2k-1,2l-1}]) &= 0,
\end{align*}
we obtain
\begin{align*}
\delta_{kj}\lambda(y_{2i,2l}) &= \delta_{il}\lambda(y_{2k,2j}), & \delta_{kj}\lambda(y_{2i-1,2l-1}) &= \delta_{il}\lambda(y_{2k-1,2j-1}), \\
\delta_{kj}\lambda(y_{2i,2l}) &= -\delta_{il}\lambda(y_{2k,2j}), & \delta_{kj}\lambda(y_{2i-1,2l-1}) &= -\delta_{il}\lambda(y_{2k-1,2j-1}).
\end{align*}
For all $1\leq i ,j\leq n$, we now have
\begin{align*}
\lambda(y_{2i-1,2j-1}) = \lambda(y_{2i,2j}) = \lambda(\theta_{2i-1,2j-1}) = \lambda(\theta_{2i,2j}) = 0.
\end{align*}
Therefore, the only character on $\mf{str}(JQ(n))$ is the trivial character $\lambda = 0$.

\subsection{The Bessel operators}\label{SS_Q_Bessel}
A technical but straightforward calculation gives
\begin{align*}
\bessel(y_{kl}) &= \sum_{i,j}y_{ij}(\pt {y_{ik}}\pt {y_{lj}}- \pt {\theta_{ik}}\pt {\theta_{lj}}) + \theta_{ij}(\pt {y_{ik}}\pt {\theta_{lj}} + \pt {\theta_{ik}}\pt {y_{lj}}), \\
\bessel(\theta_{kl}) &= \sum_{i,j}\theta_{ij}(\pt {y_{ik}}\pt {y_{lj}} + \pt {\theta_{ik}}\pt {\theta_{lj}}) + y_{ij}(\pt {y_{ik}}\pt {\theta_{lj}} - \pt {\theta_{ik}}\pt {y_{lj}}).
\end{align*}

Since we only have the trivial character on $\mf{str}(JQ(n))$, we can use the same arguments as in the periplectic case to conclude that the Bessel-Fischer product is only non-degenerate when $F_\lambda$ is isomorphic to $\ds K$ or $0$.

\subsection*{Acknowledgements}
SB is supported by a FWO postdoctoral junior fellowship from the Research Foundation Flanders (1269821N).

\bibliography{citations} 
\bibliographystyle{ieeetr}

\end{document}